\documentclass[11pt,leqno,fleqn]{article}
\usepackage{epsf,amsmath,amssymb,amscd,hyperref}
\usepackage{amsfonts}
\newcommand{\dps}{\displaystyle}

\newcommand{\GL}{\text{\rm GL}}

\newcounter{figuren}

\newcommand{\dy}[2]{%
\refstepcounter{equation}%
\LABEL{#1}%
\begin{list}{}{
\topsep 5mm
\leftmargin 18mm
\rightmargin 0cm
\itemsep 0mm
\listparindent 0mm
\parsep 0mm
\itemsep 0mm
\labelsep 0mm
\labelwidth 18mm
}%
\item[\rm (\theequation)\hfill]
#2
\end{list}%
}
\newcommand{\dyz}[1]{%
\refstepcounter{equation}%
\begin{list}{}{
\topsep 5mm
\leftmargin 18mm
\rightmargin 0cm
\itemsep 0mm
\listparindent 0mm
\parsep 0mm
\itemsep 0mm
\labelsep 0mm
\labelwidth 18mm
}%
\item[\rm (\theequation)\hfill]
#1
\end{list}%
}
\newcommand{\dyyz}[1]{\dyz{\raggedright$\dps#1$}}

\newcommand{\de}[2]{\dy{#1}{\raggedright$\displaystyle #2 $}}
\newcommand{\dez}[1]{\dyz{\raggedright$\displaystyle #1 $}}
\newcommand{\leeg}[1]{}

\newcounter{stelling}
\newcommand{\thm}[2]{\refstepcounter{stelling}\vspace{4mm}\noindent{\bf Theorem \thestelling.}\label{#1}{\it #2}}
\newcounter{bewering}
\newcommand{\prop}[2]{\refstepcounter{bewering}\vspace{4mm}\noindent{\bf Proposition \thebewering.}\label{#1}{\it #2}}
\newcounter{sectie}
\newcounter{subsectie}
\newcommand{\sect}[2]{\refstepcounter{sectie}\setcounter{subsectie}{0}
\section*{\boldmath \thesectie. #2}%
\label{#1}}
\newcommand{\sectz}[1]{\refstepcounter{sectie}\setcounter{subsectie}{0}
\section*{\boldmath \thesectie. #1}%
}
\newcommand{\subsect}[2]{\refstepcounter{subsectie}
\subsection*{\boldmath \thesectie.\thesubsectie. #2}%
\label{#1}}

\newcounter{lit}

\newcommand{\pf}{\vspace{3mm}\noindent{\bf Proof.}\ }

\newcommand{\bx}{\hspace*{\fill} \hbox{\hskip 1pt \vrule width 4pt height 8pt depth 1.5pt \hskip 1pt}

\addvspace{4mm}}

\newcommand{\rf}[1]{{\rm (\ref{#1})}}

\newcommand{\DD}{{\cal D}}

\newcommand{\GG}{{\cal G}}

\newcommand{\II}{{\cal I}}
\newcommand{\JJ}{{\cal J}}

\newcommand{\OO}{{\cal O}}

\newcommand{\NIET}[1]{}

\newcommand{\LABEL}[1]{\label{#1}}

\newcommand{\rank}{\text{\rm rank}}
\newcommand{\sgn}{\text{\rm sgn}}

\newcommand{\oC}{{\mathbb{C}}}

\newcommand{\oF}{{\mathbb{F}}}

\newcommand{\oN}{{\mathbb{N}}}

\newcommand{\oR}{{\mathbb{R}}}

\renewcommand{\Im}{{\text{Im~}}}
\newcommand{\Ker}{{\text{Ker~}}}

\begin{document}

\begin{center}
{\LARGE\bf Characterizing partition functions of the vertex model

}
\end{center}
\vspace{1mm}
\begin{center}
{\large
\hspace{10mm}
Jan Draisma\footnote{ University of Technology Eindhoven and
CWI Amsterdam
$\mbox{}^2$
CWI Amsterdam and Department of Mathematics, Leiden University
$\mbox{}^3$
Department of Computer Science, E\"otv\"os Lor\'and Tudom\'anyegyetem
Budapest
(The
European Union and the European Social Fund have provided financial support to the
project under the grant agreement no. T\'AMOP 4.2.1./B-09/KMR-2010-0003.)
$\mbox{}^4$
CWI Amsterdam
$\mbox{}^5$
CWI Amsterdam and Department of Mathematics, University of Amsterdam.
Mailing address: CWI, Science Park 123, 1098 XG Amsterdam,
The Netherlands.
Email: \url{lex@cwi.nl}
},
Dion C. Gijswijt$\mbox{}^2$,
L\'aszl\'o Lov\'asz$\mbox{}^3$,
Guus Regts$\mbox{}^4$,
and
Alexander Schrijver$\mbox{}^5$
\addtocounter{footnote}{4}

}

\end{center}

\noindent
{\small{\bf Abstract.}
We characterize which graph parameters are partition functions
of a vertex model over an algebraically closed field of
characteristic 0 (in the sense of de la Harpe and Jones [4]).
We moreover characterize when the vertex model can be taken so that
its moment matrix has finite rank.
%Basic techniques are the Nullstellensatz and the First and
%Second Fundamental Theorems of Invariant theory for the orthogonal
%group.

}

\sect{25no10a}{Survey of results}

Let $\GG$ denote the collection of all undirected graphs, two of them
being the same if they are isomorphic.
In this paper, all graphs are finite and may have loops and multiple edges.
Let $k\in\oN$ and let $\oF$ be a commutative ring.
Following de la Harpe and Jones [4], call any function $y:\oN^k\to\oF$ a
({\em $k$-color}) {\em vertex model} ({\em over $\oF$}).\footnote{ In
[8] it is called an {\em edge coloring model}.
Colors are also called {\em states}.}
The {\em partition function} of $y$ is the function $f_y:\GG\to\oF$
defined for any graph $G=(V,E)$ by
\dez{
f_y(G):=\sum_{\kappa:E\to[k]}\prod_{v\in V}y_{\kappa(\delta(v))}.
}
Here $\delta(v)$ is the set of edges incident with $v$.
Then $\kappa(\delta(v))$ is a multisubset of $[k]$,
which we identify with its incidence vector in $\oN^k$.
Moreover, we use $\oN=\{0,1,2,\ldots\}$ and for $n\in\oN$,
\dez{
[n]:=\{1,\ldots,n\}.
}

We can visualize $\kappa$ as a coloring of the edges of $G$ and
$\kappa(\delta(v))$ as the multiset of colors `seen' from $v$.
The vertex model was considered by
de la Harpe and Jones [4] as a physical model, where vertices serve
as particles, edges as interactions between particles,
and colors as states or energy levels.
It extends the Ising-Potts model.
Several graph parameters are partition functions of some vertex model,
like the number of matchings.
There are real-valued graph parameters that are partition functions
of a vertex model over $\oC$, but not over $\oR$.
(A simple one is $(-1)^{|E(G)|}$.)
%We give more background in Section \ref{25no10b}.

In this paper, we characterize which functions $f:\GG\to\oF$
are the partition function of a vertex model over $\oF$,
when $\oF$ is an algebraically closed field of characteristic 0.

To describe the characterization,
call a function $f:\GG\to\oF$ {\em multiplicative} if
$f(\emptyset)=1$ and $f(GH)=f(G)f(H)$ for all $G,H\in\GG$.
Here $GH$ denotes the disjoint union of $G$ and $H$.

Moreover,
for any graph $G=(V,E)$, any $U\subseteq V$, and any $s:U\to V$, define
\dyz{
$E_s:=\{us(u)\mid u\in U\}$ and $G_s:=(V,E\cup E_s)$
}
(adding multiple edges if $E_s$ intersects $E$).
Let $S_U$ be the group of permutations of $U$.

\thm{9no10c}{
Let $\oF$ be an algebraically closed field of characteristic $0$.
A function $f:\GG\to\oF$ is the partition function of some
$k$-color vertex model over $\oF$
if and only if $f$ is multiplicative
and for each graph $G=(V,E)$, each $U\subseteq V$ with $|U|=k+1$, and
each $s:U\to V$:
\de{10no10i}{
\sum_{\pi\in S_U}\sgn(\pi)f(G_{s\circ\pi})=0.
}
}

Let $y:\oN^k\to\oF$.
The corresponding {\em moment matrix} is
\dez{
M_y:=(y_{\alpha+\beta})_{\alpha,\beta\in\oN^k}.
}
Abusing language we say that $y$ has {\em rank} $r$ if $M_y$ has rank $r$.
For any graph $G=(V,E)$, $U\subseteq V$, and $s:U\to V$, let
$G/s$ be the graph obtained from $G_s$ by contracting all edges in $E_s$.

\thm{9no10d}{
Let $f$ be the partition function of a $k$-color
vertex model over an algebraically closed field $\oF$ of
characteristic $0$.
Then $f$ is the partition function of a $k$-color vertex model over
$\oF$ of rank at most $r$ if and only if
for each graph $G=(V,E)$, each $U\subseteq V$ with $|U|=r+1$, and
each $s:U\to V\setminus U$:
\de{10no10j}{
\sum_{\pi\in S_U}\sgn(\pi)f(G/{s\circ\pi})=0.
}
}

\medskip
It is easy to see that the conditions in Theorem \ref{9no10d} imply
those in Theorem \ref{9no10c} for $k:=r$, since for each $u\in U$ we
can add to $G$ a new vertex $u'$ and a new edge $uu'$, thus obtaining
graph $G'$.
Then \rf{10no10j} for $G'$, $U'$, and $s'(u'):=s(u)$ gives
\rf{10no10i}.

This implies that if $f$ is the partition function of a vertex
model of rank $r$, it is also the partition function of an $r$-color
vertex model of rank $r$.

It is also direct to see that in both theorems we may restrict
$s$ to injective functions.
However, in Theorem \ref{9no10c}, $s(U)$ should be allowed to
intersect $U$
(otherwise $f(G):=2^{\#\text{ of loops}}$ would satisfy the
condition for $k=1$, but is not the partition function of
some $1$-color vertex model).
Moreover, in Theorem \ref{9no10d}, $s(U)$ may not intersect $U$
(otherwise $f(G):=2^{|V(G)|}$ would not satisfy the condition for $k=r=1$,
while it is the partition function of some 1-color vertex model of
rank 1).

%\NIET{
\sect{25no10b}{Background}

In this section,
we give some background to the results described in this paper.
The definitions and results given in this section will not be used
in the remainder of this paper.

As mentioned, the vertex model roots in
mathematical physics, see de la Harpe and Jones [4].
They also gave the dual `spin model', where the roles of vertices
and edges are interchanged.
Both models are generalizations of the Ising-Potts model of statistical
mechanics.
Partition functions of spin models were characterized by
Freedman, Lov\'asz, and Schrijver [2]
and
Schrijver [7].

We describe some results of Szegedy [8,\linebreak[0]9] concerning the
vertex model that are related to, and have motivated, our results.
They require the notions of $l$-labeled graphs and $l$-fragments.

For $l\in\oN$, an
{\em $l$-labeled graph} is an undirected graph $G=(V,E)$
together with an injective `label' function $\lambda:[l]\to V$.
If $G$ and $H$ are two $l$-labeled graphs, let $GH$ be the graph
obtained from the disjoint union of $G$ and $H$ by identifying
equally labeled vertices.
(We can identify (unlabeled) graphs with $0$-labeled graphs, and then this
notation extends consistently the notation $GH$ given in
Section \ref{25no10a}.)

An {\em $l$-fragment} is an $l$-labeled graph where each
labeled vertex has degree 1.
(If you like, you may alternatively view the degree-1 vertices
as ends of `half-edges'.)
If $G$ and $H$ are $l$-fragments, the graph $G\cdot H$
is obtained from $GH$ by ignoring each of the $l$ identified points
as vertex, joining its two incident edges into one edge.
(A good way to imagine this is to see a graph as a topological $1$-complex.)
Note that it requires that we also should consider the `vertexless
loop' as possible edge of a graph, as we may create it in $G\cdot H$.

Let $\GG_l$ and $\GG'_l$ denote the collections of $l$-labeled graphs
and of $l$-fragments, respectively.
For any $f:\GG\to\oF$ and $l\in\oN$, the {\em connection matrices}
$C_{f,l}$ and $C'_{f,l}$ are the $\GG_l\times\GG_l$ and
$\GG'_l\times\GG'_l$ matrices defined by
\dyz{
$\dps C_{f,l}:=(f(GH))_{G,H\in\GG_l}$
~~~and~~~
$\dps C'_{f,l}:=(f(G\cdot H))_{G,H\in\GG'_l}$.
}
Now we can formulate Szegedy's theorem ([8]):
\dyz{
A function $f:\GG\to\oR$ is the partition function of
a vertex model over $\oR$
if and only if $f$ is multiplicative and
$C'_{f,l}$ is positive semidefinite for each $l$.
}
Note that the number of colors is equal to the $f$-value of the
vertexless loop.
The proof is based on the First Fundamental Theorem for the orthogonal
group and on the Real Nullstellensatz.

Next consider the complex case.
Szegedy [9] observed that if $y$ is a vertex model of rank $r$,
then $\rank(C_{f_y,l})\leq r^l$ for each $l$.
It made him ask whether, conversely, for each function
$f:\GG\to\oC$ with $f(\emptyset)=1$ such that there exists
a number $r$ for which $\rank(C_{f,l})\leq r^l$ for each $l$, there
exists a finite rank vertex model $y$ over $\oC$ with $f=f_y$.
The answer is negative however: the function $f$ defined by
\dez{
f(G):=
\begin{cases}
(-2)^{\#\text{ of components}}&\text{ if $G$ is $2$-regular,}\\
0&\text{ otherwise,}
\end{cases}
}
has $f(\emptyset)=1$ and can be shown to have
$\rank(C_{f,l})\leq 4^l$ for each $l$.
However, $f$ is not the partition function of a vertex model
(as it does not satisfy condition of Theorem \ref{9no10c} for any $k$).
The characterizations given in the present paper may serve as
alternatives to Szegedy's question.
%}

\sect{21no10a}{Proof of Theorem \ref{9no10c}}

We fix $k$.
Necessity of the conditions is direct.
Condition \rf{10no10i} follows from the fact that, as $|U|=k+1$, for any
$\kappa:E\cup E_s\to[k]$ there exist
distinct $u,v\in U$ with $\kappa(us(u))=\kappa(vs(v))$.
As the permutation exchanging $u$ and $v$ has negative sign, this
gives cancellation in the sum \rf{10no10i}.

To see sufficiency,
introduce a variable $y_{\alpha}$ for each $\alpha\in\oN^k$ and 
define the ring $R$ of polynomials in these variables
Define
\dez{
R:=\oF[y_{\alpha}\mid\alpha\in\oN^k].
}
There is a bijection between the variables $y_{\alpha}$
in $R$ and the monomials
$x^{\alpha}=\prod_{i\in\alpha}x_i$ in $\oF[x_1,\ldots,x_k]$.
(Note that
$x^{\alpha}x^{\beta}$ does not correspond to
$y_{\alpha}y_{\beta}$, but with $y_{\alpha+\beta}$.)
In this way, functions $y:\oN^k\to\oF$ correspond to
elements of $\oF[x_1,\ldots,x_k]^*$.

Define $p:\GG\to R$ by $p(G)(y):=f_y(G)$ for any graph $G=(V,E)$
and $y:\oN^k\to\oF$.
We must show that the polynomials $p(G)-f(G)$ have a common zero.

Let $\oF\GG$ denote the set of formal $\oF$-linear combinations
of elements of $\GG$.
The elements of $\oF\GG$ are called {\em quantum graphs}.
We can extend $f$ and $p$ linearly to $\oF\GG$.
Taking disjoint union of graphs $G$ and $H$ as product $GH$,
makes $\oF\GG$ to an algebra.
Then $f$ and $p$ are algebra homomorphisms.

The main ingredients of the proof are two basic facts about $p$
(independently of $f$): a characterization of the image and a
characterization of the kernel of $p$.

Let $I$ be the subspace of $\oF\GG$ spanned by the quantum graphs
\de{9no10a}{
\sum_{\pi\in S_U}\sgn(\pi)G_{s\circ\pi},
}
where $G=(V,E)$ is a graph, $U\subseteq V$ with $|U|=k+1$,
and $s:U\to V$.
Then
\de{17no10b}{
\Ker p=I.
}

To characterize $\Im p$, let $O_k$ be the group of orthogonal
matrices over $\oF$ of order $k$.
Observe that
$O_k$ acts on $\oF[x_1,\ldots,x_k]$, and hence on $R$,
through the bijection $y_{\alpha}\leftrightarrow x^{\alpha}$
mentioned above.
Then, as was observed by Szegedy [8],
\de{17no10c}{
\Im p=R^{O_k},
}
where as usual, $Z^{O_k}$ denotes the set of $O_k$-invariant
elements of $Z$, if $O_k$ acts on a set $Z$.

\rf{17no10c} and \rf{17no10b}
follow from the First and Second Fundamental Theorems of
Invariant Theory for $O_k$,
as we will show in Section \ref{21no10a}.\ref{10no10e}.

As $f$ is multiplicative, $f$ extends to an algebra
homomorphism $f:\oF\GG\to\oF$.
By the condition in Theorem \ref{9no10c}, $f(I)=0$.
Hence by \rf{17no10b} there exists an algebra homomorphism
$\hat f:p(\oF\GG)\to\oF$ such that $\hat f\circ p=f$.

Let $\II$ be the ideal generated by the polynomials $p(G)-f(G)$ for
graphs $G$.
Let $\rho_{O_k}$ denote the Reynolds operator.
By \rf{17no10c}, $\rho_{O_k}(\II)$ is equal to the ideal in
$p(\oF\GG)$ generated by the polynomials $p(G)-f(G)$.
(This follows essentially from the fact that if
$q\in R^{O_k}$ and $r\in R$, then
$\rho_{O_k}(qr)=q\rho_{O_k}(r)$.)
This implies, as $\hat f(p(G)-f(G))=0$,
\de{16no10a}{
\hat f(\rho_{O_k}(\II))=0,
}
hence $1\not\in\II$.

If $|\oF|$ is uncountable, the Nullstellensatz for countably
many variables (Lang [6]) yields
the existence of a common zero $y$.

To prove it for general algebraically closed fields $\oF$ of
characteristic 0, let,
for each $d\in\oN$,
$Z_d:=\{\alpha\in\oN^k\mid |\alpha|\leq d\}$
and
\dyz{
$Y_d:=\{y|Z_d\mid q(y)=\hat f(q)$ for each
$q\in\oF[y_{\alpha}\mid\alpha\in Z_d]^{O_k}\}$.
}
By the Nullstellensatz, since $|Z_d|$ is finite,
$Y_d\neq\emptyset$ for each $d$.
Note that $Y_d$ is $O_k$-stable.
This implies that $Y_d$ contains a unique $O_k$-orbit $C_d$ of minimal
(Krull) dimension
(cf.\ [5] Satz 2, page 101
or [1] 1.11 and 1.24).

Let $\pi_d$ be the projection $z\mapsto z|Z_d$ for $z$ belonging to
any $Y_{d'}$ ($d'\geq d$).
Note that if $d'\geq d$ then
$\pi_d(C_{d'})$ is an $O_k$-orbit contained in $Y_d$.
Hence
\de{22no10a}{
\dim C_d\leq\dim\pi_d(C_{d'})\leq\dim C_{d'}.
}
As $\dim C_d\leq\dim O_k$ for all $d$, there is a $d_0$ such that
for each $d\geq d_0$, $\dim C_d=\dim C_{d_0}$.
Hence we have equality throughout in \rf{22no10a} for all
$d'\geq d\geq d_0$.

By the uniqueness of the orbit of smallest dimension,
this implies that, for all $d'\geq d\geq d_0$,
$C_d=\pi_{d}(C_{d'})$.
Hence there exists $y:\oN^k\to\oF$ such that
$y|Z_d\in C_d$ for each $d\geq d_0$.
This $y$ is as required.

\subsect{10no10e}{Applying the Fundamental Theorems for $O_k$}

Let $n\in\oN$, and
let $\GG_n$ be the collection of graphs with $n$ vertices.
Let $S\oF^{n\times n}$ be the set of symmetric matrices in
$\oF^{n\times n}$.
For any linear space $X$, let $\OO(X)$ denote the space of regular
functions on $X$ (the algebra generated by the linear functions
on $X$).
Then $\OO(S\oF^{n\times n})$ is spanned by the monomials
$\prod_{ij\in E}x_{i,j}$ in the variables $x_{i,j}$,
where $([n],E)$ is a graph.
Here $x_{i,j}=x_{j,i}$ are the standard coordinate functions on
$S\oF^{n\times n}$, taking $ij$ as unordered pair.

We define linear functions $\mu$, $\sigma$, and $\tau$ so that
the following diagram commutes:
\de{30no10a}{
\begin{CD}
\oF\GG_n @>p>> R_n\\
@AA\mu A       @AA\sigma A\\
\OO(S\oF^{n\times n}) @>\tau>>  \OO(\oF^{k\times n})
\end{CD},
}
Here $\oF\GG_n$ is the linear space of formal linear combinations
of elements of $\GG_n$,
and $R_n$ is the set of homogeneous polynomials in $R$ of degree $n$.

Define $\mu$ by
\dy{21no10c}{
$\mu(\prod_{ij\in E}x_{i,j}):=G$
}
for any graph $G=([n],E)$.
Define $\sigma$ by
\dy{20no10b}{
$\dps\sigma(\prod_{j=1}^n\prod_{i=1}^kz_{i,j}^{\alpha(i,j)})
:=\prod_{j=1}^ny_{\alpha_j}$
}
for $\alpha\in\oN^{k\times n}$, where $z_{i,j}$ are the standard
coordinate functions on $\oF^{k\times n}$ and
where $\alpha_j=(\alpha(1,j),\ldots,\alpha(k,j))\in\oN^k$.
Then $\sigma$ is $O_k$-equivariant.

Finally, define $\tau$ by
\dy{8de10a}{
$\dps\tau(q)(z):=q(z^Tz)$
}
for $q\in\OO(S\oF^{n\times n})$ and $z\in\oF^{k\times n}$.

Now \rf{30no10a} commutes; in other words,
\de{20no10a}{
p\circ\mu=\sigma\circ\tau.
}
To prove it, consider any monomial $q:=\prod_{ij\in E}x_{i,j}$ in
$\OO(S\oF^{n\times n})$, where $G=([n],E)$ is a graph.
Then for $z\in\oF^{k\times n}$,
\dyyz{
\tau(q)(z)=q(z^Tz)=\prod_{ij\in E}\sum_{h=1}^kz_{h,i}z_{h,j}=
\sum_{\kappa:E\to[k]}\prod_{i\in[n]}\prod_{e\in\delta(i)}z_{\kappa(e),i}.
}
So, by definition \rf{20no10b} of $\sigma$ and \rf{21no10c} of $\mu$,
\dyyz{
\sigma(\tau(q))=
\sum_{\kappa:E\to[k]}\prod_{i\in[n]}y_{\kappa(\delta(i))}=p(G)
=p(\mu(q)).
}
This proves \rf{20no10a}.

Note that $p$ and $\tau$ are algebra homomorphisms, but
$\mu$ and $\sigma$ generally are not.
The latter two functions are surjective,
and their restrictions to the $S_n$-invariant part of their 
respective domains are bijective.

The First Fundamental Theorem (FFT) for $O_k$
%(rrr[De Concini,Procesi,1976] Theorem 5.6,)
(cf.\ [3] Theorem 5.2.2) says that
$\Im \tau=(\OO(\oF^{k\times n}))^{O_k}$.
Hence, as $\mu$ and $\sigma$ are surjective, and as $\sigma$ is
$O_k$-equivariant, $\Im p=R_n^{O_k}$.
This implies \rf{17no10c}.

The Second Fundamental Theorem (SFT) for $O_k$
%(rrs[De Concini,Procesi,1976] Theorem 5.7,)
(cf.\ [3] Theorem 12.2.14) says that
$\Ker\tau=K$,
where $K$ is the ideal in $\OO(S\oF^{n\times n})$ generated by the
$(k+1)\times(k+1)$ minors of $S\oF^{n\times n}$.

This implies $\Ker p=I$.
Indeed, $I\subseteq\Ker p$ follows from the necessity of the
conditions of Theorem \ref{9no10c}.
To see the reverse inclusion, let $\gamma\in\oF\GG$ with $p(\gamma)=0$.
We can assume $\gamma\in\oF\GG_n$.
Then $\gamma=\mu(q)$ for some $q\in(\OO(S\oF^{n\times n}))^{S_n}$.
Hence $\sigma(\tau(q))=p(\mu(q))=p(\gamma)=0$.
As $\tau(q)$ is $S_n$-invariant, this implies $\tau(q)=0$.
So $q\in K$, hence $\gamma=\mu(q)\in\mu(K)\subseteq I$.
This gives \rf{17no10b}.

\sectz{Proof of Theorem \ref{9no10d}}

Necessity can be seen as follows.
Choose $y:\oN^k\to\oF$ with $\rank(M_y)\leq r$
and choose $\kappa:E\to [k]$, $U\subseteq V$ with $|U|=r+1$,
and $s:U\to V\setminus U$.
Then
\dyyz{
\sum_{\pi\in S_U}\sgn(\pi)f_y(G/s\circ\pi)
=
\sum_{\kappa:E\to[k]}
\sum_{\pi\in S_U}\sgn(\pi)
\prod_{u\in U}
y_{\kappa(\delta(u)\cup\delta(s(\pi(u))))}
\cdot
\hspace{-6mm}
\prod_{v\in V\setminus(U\cup s(U))}
\hspace{-6mm}
y_{\kappa(\delta(v))}
=
\sum_{\kappa:E\to[k]}
\det(y_{\kappa(\delta(u)\cup\delta(s(v)))})_{u,v\in U}
\hspace{-6mm}
\prod_{v\in V\setminus(U\cup s(U))}
\hspace{-6mm}
y_{\kappa(\delta(v))}
=
0.
}

To see sufficiency,
let $J$ be the ideal in $\oF\GG$ spanned by the quantum graphs
\de{9no10e}{
\sum_{\pi\in S_U}\sgn(\pi)G/{s\circ\pi},
}
where $G=(V,E)$ is a graph, $U\subseteq V$ with $|U|=r+1$, and
$s:U\to V\setminus U$.
%Note that $I\subseteq J$ (by the reduction from \rf{10no10i} to \rf{10no10j}
%given in Section \ref{25no10a}).
Let $\JJ$ be the ideal in $R$ generated by the polynomials
$\det N$ where $N$ is a $(r+1)\times(r+1)$ submatrix of $M_y$.

\prop{2no07b}{
$\rho_{O_k}(\JJ)\subseteq p(J)$.
}

\pf
It suffices to show that for any $(r+1)\times(r+1)$ submatrix
$N$ of $M_y$ and any monomial $a$ in $R$, $\rho_{O_k}(a\det N)$
belongs to $p(J)$.
Let $a$ have degree $d$, and let $n:=2(r+1)+d$.
Let $U:=[r+1]$ and let $s:U\to [n]\setminus U$ be defined by
$s(i):=r+1+i$ for $i\in[r+1]$.

We use the framework of Section \ref{21no10a}.\ref{10no10e}, with
$\tau$ as in \rf{8de10a}.
For each $\pi\in S_{r+1}$ we define linear function
$\mu_{\pi}$ and $\sigma_{\pi}$
so that the following diagram commutes:
\de{8de10b}{
\begin{CD}
\oF\GG_m @>p>> R_m\\
@AA\mu_{\pi} A       @AA\sigma_{\pi} A\\
\OO(S\oF^{n\times n}) @>\tau>>  \OO(\oF^{k\times n})
\end{CD},
}
where $m:=r+1+d$.

The function $\mu_{\pi}$ is defined by
\dyz{
$\dps\mu_{\pi}(\prod_{ij\in E}x_{i,j}):=G/s\circ\pi$
}
for any graph $G=([n],E)$.
It implies that for each $q\in \OO(S\oF^{n\times n})$,
\de{22no10g}{
\sum_{\pi\in S_{r+1}}\sgn(\pi)\mu_{\pi}(q)\in J,
}
by definition of $J$.

Next $\sigma_{\pi}$ is defined by
\dyz{
$\dps\sigma_{\pi}(\prod_{j=1}^n\prod_{i=1}^kz_{i,j}^{\alpha_{i,j}})
:=
\prod_{j=1}^{r+1}y_{\alpha_j+\alpha_{r+1+\pi(j)}}\cdot
\prod_{j=2r+3}^ny_{\alpha_i}$
}
for any $\alpha\in\oN^{k\times n}$.
So
\de{22no10h}{
a\det N=\sum_{\pi\in S_{r+1}}\sgn(\pi)\sigma_{\pi}(u)
}
for some monomial $u\in\OO(\oF^{k\times n})$.
Note that $\sigma_{\pi}$ is $O_k$-equivariant.

Now one directly checks that diagram \rf{8de10b} commutes, that is,
\de{22no10f}{
p\circ\mu_{\pi}=\sigma_{\pi}\circ\tau.
}
By the FFT, $\rho_{O_k}(u)=\tau(q)$ for some
$q\in \OO(S\oF^{n\times n})$.
Hence $\sigma_{\pi}(\rho_{O_k}(u))=
\sigma_{\pi}(\tau(q))=p(\mu_{\pi}(q))$.
Therefore, using \rf{22no10h} and \rf{22no10g},
\dez{
\rho_{O_k}(a\det N)
=
\sum_{\pi\in S_{r+1}}\sgn(\pi)\sigma_{\pi}(\rho_{O_k}(u))
=
\sum_{\pi\in S_{r+1}}\sgn(\pi)p(\mu_{\pi}(q))\in p(J),
}
as required.
\bx

(In fact equality holds in this proposition, but we do not
need it.)

Since $f$ is the partition function of a $k$-color vertex model,
there exists $\hat f:R\to\oF$ with $\hat f\circ p=f$.
If the condition in Theorem \ref{9no10d} is satisfied,
then $f(J)=0$, and hence with Proposition \ref{2no07b}
\dez{
\hat f(\rho_{O_k}(\JJ))\subseteq\hat f(p(J))=f(J)=0.
}
With \rf{16no10a} this implies that $1\not\in\II+\JJ$,
where $\II$ again is the ideal generated by the polynomials
$p(G)-f(G)$ ($G\in\GG$).
Hence $\II+\JJ$ has a common zero, as required.

\sectz{Analogues for directed graphs}

Similar results hold for directed graphs, with similar proofs,
now by applying the FFT and SFT for $\GL(k,\oF)$.
The corresponding models were also considered by de la Harpe and Jones [4].
We state the results.

Let $\DD$ denote the collection of all directed graphs, two of them being
the same if they are isomorphic.
Directed graphs are finite and may have loops and multiple edges.

The {\em directed partition function} of a $2k$-color vertex model $y$
is the function $f_y:\DD\to\oF$
defined for any directed graph $G=(V,E)$ by
\dez{
f_y(G):=\sum_{\kappa:E\to[k]}\prod_{v\in V}y_{\kappa(\delta^-(v)),\kappa(\delta^+(v))}.
}
Here $\delta^-(v)$ and $\delta^+(v)$ denote the sets of arcs
entering $v$ and leaving $v$, respectively.
Moreover, $\kappa(\delta^-(v)),\kappa(\delta^+(v))$ stands for the
concatenation of the vectors
$\kappa(\delta^-(v))$ and $\kappa(\delta^+(v))$ in $\oN^k$, so as
to obtain a vector in $\oN^{2k}$.

Call a function $f:\DD\to\oF$ {\em multiplicative} if
$f(\emptyset)=1$ and $f(GH)=f(G)f(H)$ for all $G,H\in\DD$.
Again, $GH$ denotes the disjoint union of $G$ and $H$.

Moreover,
for any directed
graph $G=(V,E)$, any $U\subseteq V$, and any $s:U\to V$, define
\dyz{
$A_s:=\{(u,s(u))\mid u\in U\}$ and $G_s:=(V,E\cup A_s)$.
}

\thm{9no10cx}{
Let $\oF$ be an algebraically closed field of characteristic $0$.
A function $f:\DD\to\oF$ is the directed partition function of some
$2k$-color vertex model over $\oF$
if and only if $f$ is multiplicative
and for each directed graph $G=(V,E)$, each $U\subseteq V$ with $|U|=k+1$,
and each $s:U\to V$:
\de{10no10ix}{
\sum_{\pi\in S_U}\sgn(\pi)f(G_{s\circ\pi})=0.
}
}

For any directed graph $G=(V,E)$, $U\subseteq V$, and $s:U\to V$, let
$G/s$ be the directed graph obtained from $G_s$ by contracting all edges in $A_s$.

\thm{9no10dx}{
Let $f$ be the directed partition function of a $2k$-color
vertex model over an algebraically closed field $\oF$ of
characteristic $0$.
Then $f$ is the directed partition function of a $2k$-color vertex model
over $\oF$ of rank at most $r$ if and only if
for each directed graph $G=(V,E)$, each $U\subseteq V$ with $|U|=r+1$, and
each $s:U\to V\setminus U$:
\de{10no10jx}{
\sum_{\pi\in S_U}\sgn(\pi)f(G/{s\circ\pi})=0.
}
}

%\sectz{Some further observations}
%
%Note that by \rf{17no10b} and \rf{17no10c}:
%\dez{
%R^{O_k}=p(\oF\GG)\cong\oF\GG/I.
%}
%
\section*{References}\label{REF}
{\small
\begin{itemize}{}{
\setlength{\labelwidth}{4mm}
\setlength{\parsep}{0mm}
\setlength{\itemsep}{1mm}
\setlength{\leftmargin}{5mm}
\setlength{\labelsep}{1mm}
}
\item[\mbox{\rm[1]}] M. Brion, 
Introduction to actions of algebraic groups,
{\em Les cours du C.I.R.M.} 1 (2010) 1-22.

\item[\mbox{\rm[2]}] M.H. Freedman, L. Lov\'asz, A. Schrijver, 
Reflection positivity, rank connectivity, and homomorphisms of graphs,
{\em Journal of the American Mathematical Society} 20 (2007) 37--51.

\item[\mbox{\rm[3]}] R. Goodman, N.R. Wallach, 
{\em Symmetry, Representations, and Invariants},
Springer, Dordrecht, 2009.

\item[\mbox{\rm[4]}] P. de la Harpe, V.F.R. Jones, 
Graph invariants related to statistical mechanical models:
examples and problems,
{\em Journal of Combinatorial Theory, Series B} 57 (1993) 207--227.

\item[\mbox{\rm[5]}] H. Kraft, 
{\em Geometrische Methoden in der Invariantentheorie},
Vieweg, Braunschweig, 1984.

\item[\mbox{\rm[6]}] S. Lang, 
Hilbert's Nullstellensatz in infinite-dimensional space,
{\em Proceedings of the American Mathematical Society} 3 (1952) 407--410.

\item[\mbox{\rm[7]}] A. Schrijver, 
Graph invariants in the spin model,
{\em Journal of Combinatorial Theory, Series B} 99 (2009) 502--511.  

\item[\mbox{\rm[8]}] B. Szegedy, 
Edge coloring models and reflection positivity,
{\em Journal of the American Mathematical Society}
20 (2007) 969--988.

\item[\mbox{\rm[9]}] B. Szegedy, 
Edge coloring models as singular vertex coloring models,
in: {\em Fete of Combinatorics and Computer Science}
(G.O.H. Katona, A. Schrijver, T.Sz\H onyi, editors),
Springer, Heidelberg and J\'anos Bolyai Mathematical Society,
Budapest, 2010, pp. 327--336.

\end{itemize}
}

\end{document}